\documentclass[12pt]{article}
\usepackage[margin=1.0in]{geometry}
\usepackage{amsmath, amssymb, amsthm, upref, rotate, epsf,array,  scalerel, enumerate, upref}

\usepackage{graphicx} 
\usepackage[utf8]{inputenc}
\usepackage{amsfonts}
\usepackage{amsmath}
\usepackage{amssymb}
\usepackage{setspace}
\usepackage{CJK}
\usepackage{moresize}
\usepackage[english]{babel}
\usepackage{xcolor}

\usepackage{mathtools}
\usepackage[utf8]{inputenc}
\usepackage[T1]{fontenc}

\theoremstyle{definition}

\newcounter{example}

\usepackage{graphicx, color}
\usepackage[sort, numbers]{natbib}
\usepackage{parskip}
\usepackage{abstract}

\date{} 

\title{\textbf{Reconstructions of piece-wise continuous and discrete functions using moments}}
\author{Robert M. Mnatsakanov{\footnote{School of Mathematical and Data Sciences, West Virginia University, Morgantown, WV 26506, USA}}, Rafik H. Aramyan{\footnote{Institute of Mathematics, National Academy of Sciences, Yerevan,   Armenia}}, and Farhad Jafari{\footnote{Department of Radiology, University of Minnesota, Minneapolis, MN 55455, USA}}
}

\begin{document}

\maketitle

%
%
%
%
%

\centerline{\bf Abstract}

\medskip

\centerline{
\parbox[b]{32pc}
{\small The problem of recovering  a moment-determinate multivariate function $f$ via its moment sequence is studied.  Under mild conditions on $f$, the point-wise and $L_1$-rates of convergence for the proposed constructions  are established. The cases where $f$ is the indicator function of a set,  and represents a discrete  probability mass function are also investigated. Calculations of the approximants and simulation studies are conducted to graphically  illustrate  the behavior of the approximations in several simple examples. Analytical and simulated errors of proposed  approximations are recorded in Tables 1-3.}
}



\noindent {{\small {\it Keywords}:   Moment-recovered approximation, rate of approximation. }


\section{Introduction and Preliminaries}

The aim of the present article is to recover a multivariate function $f$ (probability density function) by means of its moments, to describe the rate of convergence of these approximations  in various norms, and to investigate these estimates in the presence of discontinuities and in discrete settings. Proposed approximations can be used in the case of any finite $d$-dimensional model with $d \geq 2$, but for simplicity of notation we only consider the case $d=2$.

Recovery of an indicator function over some region on a plane or in space is  related to the image reconstruction problem. If the image is a polygon, its shape can be reconstructed from the knowledge of few of its moments (see Cuyt {\it et al.} \cite{Cuyt2005} and the references therein). In the latter work the so-called Pad$\acute{e}$ approximation  is used.

Recently and closely related, Henrion, Kora and Lasserre \cite{henrion2023} have used an argmin minimizer method to approximate multivariate compactly supported functions (with possible discontinuities) by polynomials $p$. Using sums of squares in the construction of $ p$, they arrive at a highly structured convex semidefinite program that provides an exact solution (with low degree polynomials), when $ f $ is a piecewise polynomial. They also mention that while exact recovery by a polynomial argmin cannot be guaranteed in general, their numerical scheme allows obtaining $ p $ when $ f $ is known by finitely many of its values. Moment methods are deeply embedded in this work and it promises to be of great significance. 

Reconstruction of a planar convex domain  from noisy observations of its moments has been studied by Goldenshluger and Spokoiny \cite{goldenshluger2004}. Using the property of one-to-one correspondence between the boundaries of convex domain and the envelopes of their support lines, they suggested
recovery of  the boundary of an image $G$ by  estimating  the support functions of $G$ from the data. In their paper, it is assumed that the Gaussian noisy  observations of the moments   of $G$ are available. The rate of convergence is derived as well.

In the current  paper we suggest the moment-recovered (MR) approximations of a bivariate function $f$  (both continuous and  discontinuous) provided that the sequence of its moments up to some finite order is known.

Our approach is general and can  be  applied for recovering any moment-determinate ($M$-determinate) function with compact support  $K \subseteq [0, T]^2, T < \infty,$ via its moments, in particular when $f$ is piece-wise constant, say, the indicator function of a convex set $G \subseteq [0, T]^2$ (see subsection 3.2).  For simplicity, we will assume that $T=1$.

To recover $f: [0, 1]^2 \to  \mathbb{R}_+$, let us use the approximation studied in   Mnatsakanov \cite{mnatsakanov2011}  and Choi et al. \cite{Choi2020} . To describe our construction consider the ordinary (geometric) moments of $f$:
\begin{align}
m_f (k, j)=  \int_0^1 \int_0^1 \,x^k \,  y^j f(x, y) \,dx \,dy := (\mathcal{B} f)  (k, j) \;\;\; k, j \in \mathbb{N}=\{0, 1, \dots \}.
\end{align}
It is known that if $ f $ is $M$-determinate, the moment sequence  ${m_f} = \{m_f ({k, j}),  k, j \in \mathbb{N}\}$ determines $f$ uniquely, (see characterizations of $ M$-determinate functions in  Shohat and Tamarkin \cite{Shohat} and Akhiezer  \cite{Akheizer2020}). To approximate the inverse transform of $\mathcal{B}$,  we introduce the sequence of operators $\{\mathcal{B}_{a}^{-1}, a=(\alpha, \alpha^{\prime})\in  \mathbb{N}_+ \times  \mathbb{N}_+ \}$:
\begin{align}\label{Eq:2}
 (\mathcal{B}_{a}^{-1} { m_f}) ( {\bf x})&:= \frac{ \Gamma(\alpha  + 2) \,
\Gamma(\alpha^{\prime}  + 2)}{ \Gamma ([\alpha x] + 1) \,  \Gamma
([\alpha^{\prime}  y] + 1)}  \\\notag
&\times\, \sum_{k=0}^{\alpha  - [\alpha x] } \,
\, \sum_{j=0}^{\alpha^{\prime}   - [\alpha^{\prime}  y] }
\frac{(-{1})^{k+j} \, m_f ({k+[\alpha x], j+[\alpha^{\prime} y]})}{k!
\, j! \, (\alpha  - [\alpha x]  - k)! \, (\alpha^{\prime} -
[\alpha^{\prime}  y]  - j)!}, \  {\bf x}=(x, y) \in [0,1]^2,
\end{align}
and construct  the MR-approximation  of $f$ defined by ${f}_{ a} :=\mathcal{B}_{a}^{-1} { m_f}$.  Here  ${m_f} = \{m_f (k, j), \\
 (k, j) \in \mathbb{N} _{a}\}$ with
$\mathbb{N} _{ a}=  \mathbb{N}_{\alpha} \times  \mathbb{N}_{\alpha^{\prime}}$ and  $\mathbb{N}_{\alpha}=  \{0, 1, \dots, \alpha \}$, while ${a}=(\alpha, \alpha^{\prime})$ with $\alpha$  and $\alpha^{\prime}  \to \infty$.  
In the sequel,  by ``$\longrightarrow_{L_1}$"  we denote the convergence in $L_1$-norm. 

\section{Recovery of continuous functions}

If the target  function $f: [0, 1]^2 \to \mathbb {R}$ is  $M$-determinate, the approximation ${f}_{a}$ specified by (2) is  analytically expressed in terms of the moments of $f$ up to $\alpha + \alpha^{\prime}$ order,   and approaches $f$  in the sup-norm as $\alpha, \alpha^{\prime} \to\infty$.
In particular, when function $f$ is sufficiently smooth  (has continuous partial  derivatives up to the second order, $f^{}_{kl}, k,l \in\{0,1, 2\}$, the rate of convergence of order $\frac{1}{\alpha}$ is derived   in \cite{Choi2020}
 as $\alpha=\alpha^{\prime}\to \infty$.  Namely,  the following statement is true.

 \noindent
 {\bf Proposition 1.} {\it Let $ f \in C^2( [0,1]^2) $.  If $\alpha=\alpha^{\prime}\to\infty$, then  
\begin{align}\label{eq:4-2}
||f_a - f||_{\infty}
& \leq  \frac{C}{\alpha+2} + o\Big (\frac{1}{\alpha}\Big ),
\end{align}
 where $||\cdot||_{\infty}$ is the sup-norm and $C= 2 \Big (||f^{}_{10}|| +  ||f^{}_{01}||\Big )  + \frac{1}{2}\, \Big ( ||f^{}_{20}|| + || f^{}_{11}|| +  ||f^{}_{02}|| \Big )$}.

\noindent
\textsc {\bf  Example 1}. Assume that $\mathcal{P}_{p+q}$ denotes the family of polynomials up to  order $p+ q \geq 0$, defined on $[0, 1]^2$ and specified by  the coefficients $a_{k,l}$.  The rate of approximation for functions $f \in \mathcal{P}_{p+q}$ were derived in  Mnatsakanov \cite{mnatsakanov2011}.  In particular, applying Corollary 2 of Theorem 2, one can recover   $f(x,y) = 3( x^2
+y^2) / 2, \, (x, y) \in [0, 1]^2$  by means of   ${f}_{a}$. In this case $p=q=2$ with only nonzero coefficients $a_{20}=a_{02}=3/2$. Note that the sequence of moments $m_f$ used in (2) has a  simple form,  $m_f = \{m_f ({k,j})=3 [(k+3)^{-1}\,
(j+1)^{-1} + (k+1)^{-1}\, (j+3)^{-1}] / 2,  (k, j) \in \mathbb{N}_{a}\}$.  Table 1  provides several values of  errors both in the sup-norm  $\Delta_{ a}=
\mid \mid f_{a}-f \mid \mid_{\infty}$ and in $L_1$-norm, $d_1 (f_a, f)$,  for the approximants ${f}_{ a}$ calculated by means of $m_f$, when the values $\alpha=\alpha^{\prime} \in \{10, 15, 25, 32,  50\}$.
Application of  Corollary 2$(i)$ of Theorem 2 in \cite{mnatsakanov2011}, provides the sup-norm of the difference
\begin{align*}
f_{ a} (x, y) - f(x, y) = \tfrac {3}{2} (b_{20} (x) - x^2) +  \tfrac{3}{2} (b_{02} (y) - y^2) \;\;\; {\rm for}\;\;\;  (x, y) \in [0, 1]^2.
\end{align*}
Here
\begin{align*}
b_{20} (x)=\frac{([\alpha x] +1)([\alpha x] +2)}{(\alpha+2)(\alpha +3)} \;\;\; {\rm and}\;\;\; b_{02} (y)=\frac{([\alpha^{\prime}y] +1)([\alpha^{\prime} y] +2)}{(\alpha^{\prime}+2)(\alpha^{\prime} +3)}.
\end{align*}
It is easily seen that the uniform error is achieved at $(x, y)=(1, 1)$, and  is exactly equal to $6/(\alpha +3)$, i.e., the exact analytical error $6/(\alpha +3)$ almost coincides  with   the actual error  $\Delta_{ a}$ recorded in Table 1. 

\vspace{0.2cm}

\begin{table}
\centering
 \caption{The  $d_1 (f_a, f)$ approximation errors, the uniform $\Delta_{ a}$-errors and simulated uniform $\widehat \Delta_{a}$-errors}
 \vspace{0.10in}
{\begin{tabular}{|c|c|c|c|c|c|c|}
  \hline
  $\alpha=\alpha^{\prime}$ & 10 & 15 & 25 & 32 & 50 \\
  \hline
  $d_1 (f_a, f)\times 10^4$ & 0.0382 & 0.0193 & 0.0076 & 0.0047 & 0.0019   \\
 $\Delta_{a}=\mid \mid f_{a}-f \mid \mid_{\infty}$              & 0.4615 & 0.3333 & 0.2143 & 0.1714 & 0.1132   \\
 $\hat \Delta_a$ with $n=10^4$ &0.4734 & 0.3824 & 0.3780 & 0.4083 & 0.5878   \\
 $\hat \Delta_a$ with $n=10^5$ &0.4582 & 0.3445 & 0.2457 & 0.2219 & 0.2264  \\
\hline
\end{tabular}}
\label{table1}
\end{table}

To demonstrate the performance of our construction in the bivariate model,  let us consider  the empirical moments   $\hat{m_f} =\{\hat m_f ({k, j}), (k, j) \in \mathbb{N}_{a}\}$ with
\begin{align*}
\hat m_f({k, j})=\frac{1}{n} \sum_{i=1}^n X_i^k \, Y_i^j \, .
\end{align*}
Applying (2) with $\hat m_f ({k, j})$ instead  of $m_f ({k, j})$ yields the MR-estimator $\hat f_{a}$ of $f$. We
simulated the data $(X_i, Y_i), i=1, \dots , n$, from $f(x, y)= 3\, (x^2 + y^2) / 2$ with the sample sizes $n=10^4$ and $10^5$.
In Figure 1(a), (b) and (c), we plotted the graphs of $f$, $f_{a}$, and $\hat f_{a}$ when $\alpha=\alpha^\prime= 25$ and $n=10^4$. In addition, we repeated the  simulations 200 times ($N=200$) when $n=10^4$ and $n=10^5$, and recorded (see Table 1) the values of simulated uniform errors
\begin{align*}
\widehat\Delta_{ a} = \frac{1}{N}\, \sum_{j=1}^N\, \mid \mid \hat f^{j}_{ a}-f \mid \mid_{\infty}\, ,
\end{align*}
for  $\alpha=\alpha^{\prime} \in \{10, 15,  25, 32,  50\}$. Here $\hat f^{j}_{a}$ denotes the empirical counterpart of  MR-approxima-
tion derived on the $j$-th replication.
From Figure 1(b) we see that the performance of $f_{a}$ is quite  good for $\alpha=25$  except at ${\bf x}=(1, 1)$.
This explains the presence  of large errors in sup-norm recorded in Table 1 when $\alpha$ and $\alpha^\prime$ are small.
Also from Table 1 we conclude  that when $n=10^4$ and $n=10^5$ the optimal values for $\alpha=\alpha^{\prime}$ are 25 and 32, respectively, i.e., the optimal value of  $\alpha=\alpha^{\prime}$  increases but very slowly as $n$ does.

\begin{figure}
\begin{center}
\begin{tabular}{ccc}
 \includegraphics[width=0.3\textwidth ]{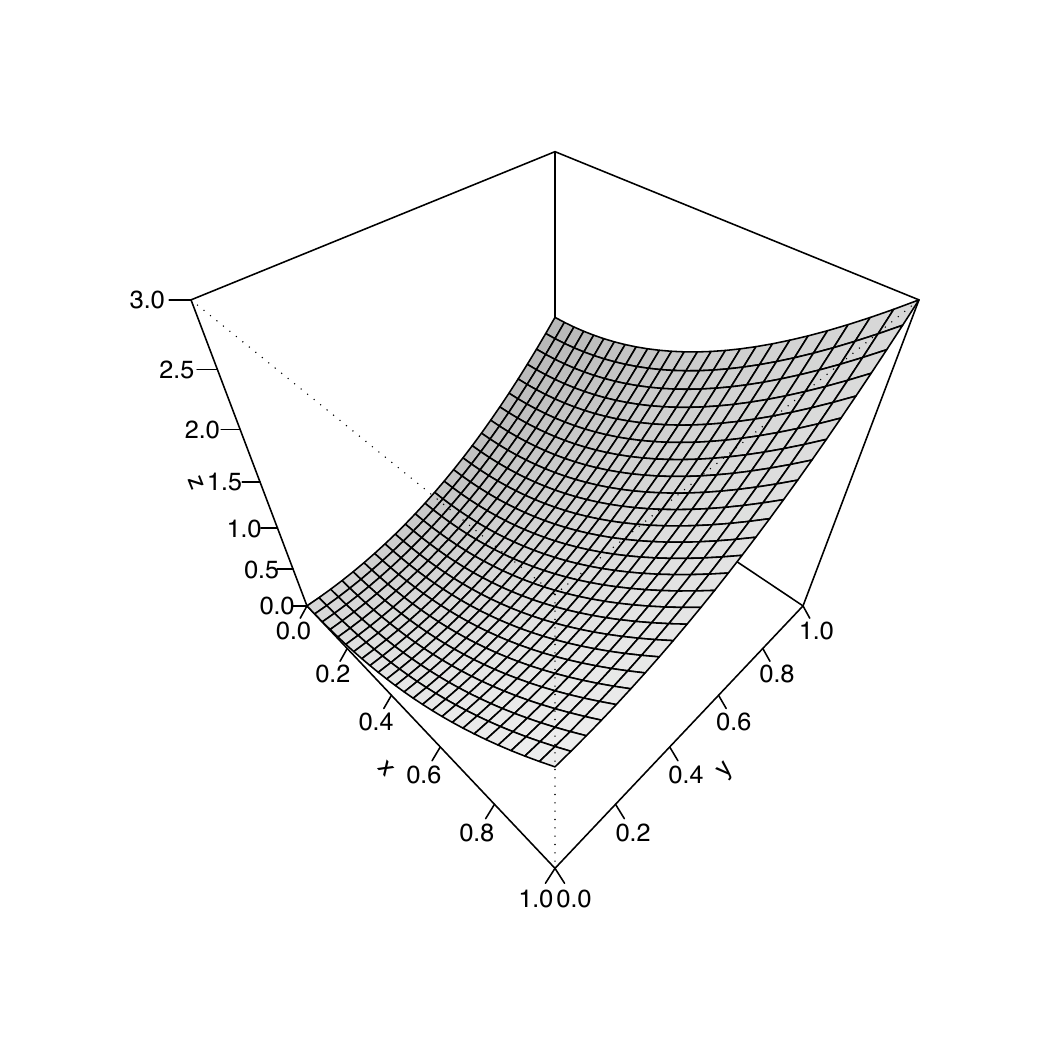} &
 \includegraphics[width=0.3\textwidth ]{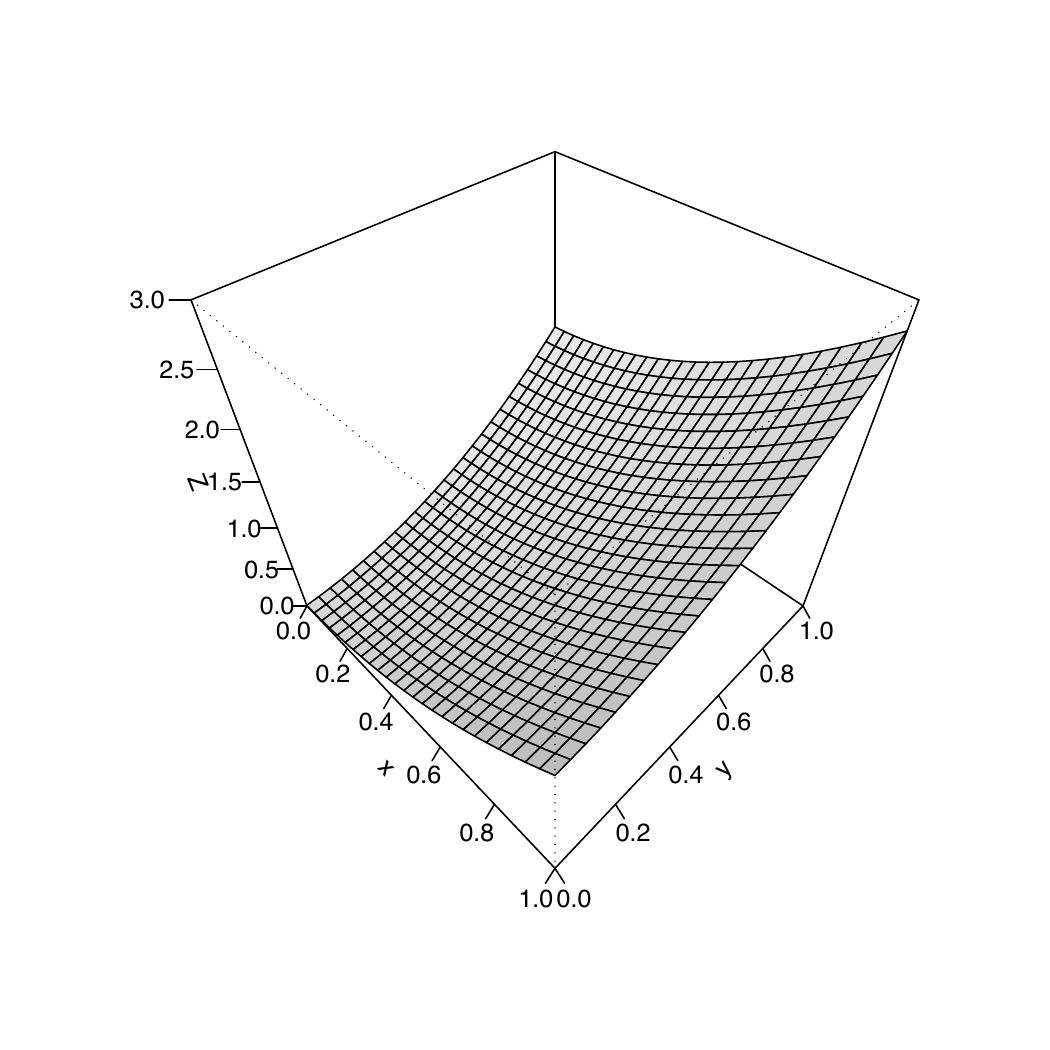}&
  \includegraphics[width=0.3\textwidth ]{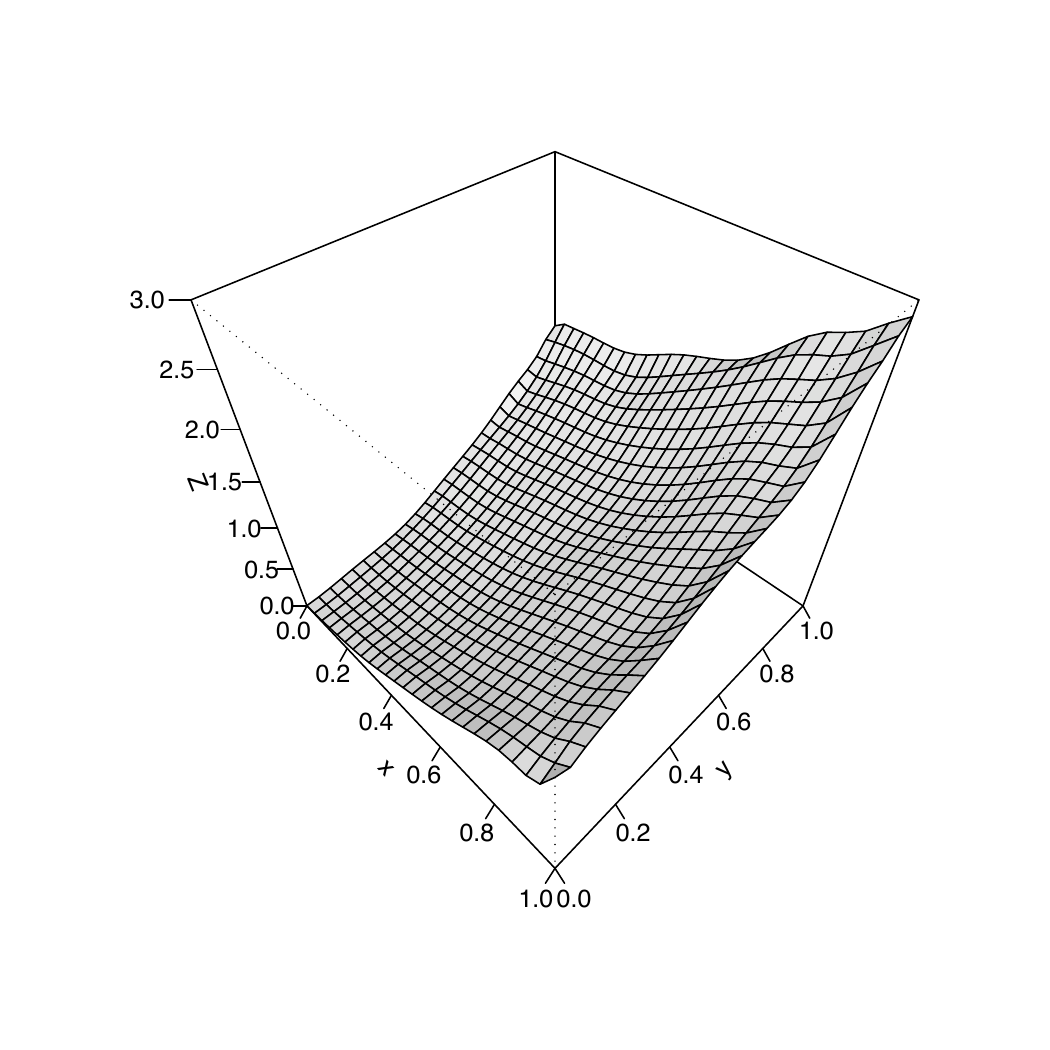}  \\
(a)  & (b) & (c)
\end{tabular}
\caption{(a) The graph  of $f(x, y)= 3 (x^2+y^2)/2$; (b) the graph
of approximant $f_{a}$ with $a=(\alpha, \alpha^\prime )$ and $\alpha=\alpha^\prime =25$; and (c)
the graph of  $\hat f_{ a}$ with $\alpha=\alpha^\prime =25$ and
$n=10^4$} \label{Figure 1}
\end{center}
\end{figure}

\noindent
{\bf Remark 1.}
If $f$ is observed directly via the sample $(X_i, Y_i), i=1, \dots , n$,  the MR-estimator,  $\hat f_{a}$,  can be rewritten  in  the form of asymmetric beta kernel density estimate (cf.  with Mnatsakanov {\it et al.} \cite{mnatsakanov2022}):
\begin{align*}
\hat f_{ a} (x, y)=\int_{[0,1]^2} \,   \beta_\alpha (t, x) \, \beta_{\alpha^\prime} (s,  y) \, d \hat F(t, s).
\end{align*}
Here $\hat F$ is the empirical cdf of $(X_i, Y_i)$'s,  while $\beta_{\alpha} (\cdot, x), \alpha\in \mathbb {N}_+,$ denotes the sequence of  beta  kernel  density functions having the shape parameters $[\alpha x] +1$, and $\alpha- [\alpha x] +1$, respectively. The reader is referred  to Bouezmarni and Rolin \cite{bouezmarni2003},  and Chen \cite{chen1999}, where the  asymmetric beta kernels density estimation  is introduced and  its consistency is studied.  Since  $E \hat f_{a}=f_{ a}$, the statements of Theorem 5.2 from \cite{Choi2020}  provides the rates of convergence for the bias term of $\hat f_{ a} $ when the underlying pdf $f$ is sufficiently smooth. The investigation of the trade off between the statistical error and approximation error of $ \hat{f}_{a} $ is postponed. 


\section{Image Reconstruction}
When the function $f({\bf x})$ represents the intensity of the image at a spatial position ${\bf x}=(x, y)$, the problem of recovering  $f$ is known as the image reconstruction problem. 
Since images have finite extensions  it is common to assume that the support of $f$ is compact, say, $[0, 1]^2$. Also, in practice, $ G$ is often a connected regular (often convex, but not necessarily here) set in the plane. To recover an image  $G$ one can use the geometric moments  (cf. with (1)) of $G$:
\begin{equation*}
m_G ({k, j}) = \int\int_G t^k s^j \, dt\, ds
\vspace{0.3cm}
\end{equation*}
instead of $m_f ({k, j})$ in (2). Consider the following  example.

\begin{figure}
\begin{center}
\begin{tabular}{cc}
\includegraphics[width=0.23\textwidth ]{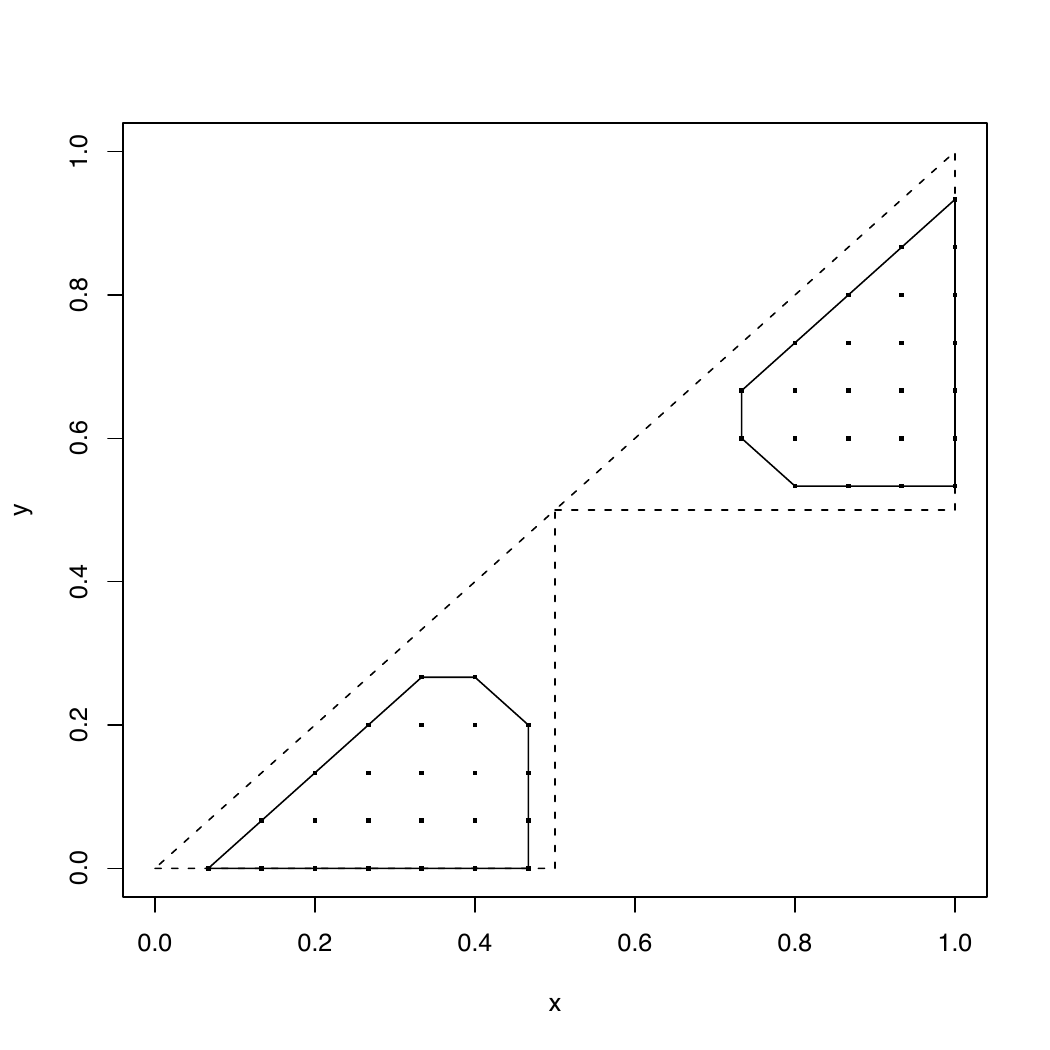}
\includegraphics[width=0.23\textwidth ]{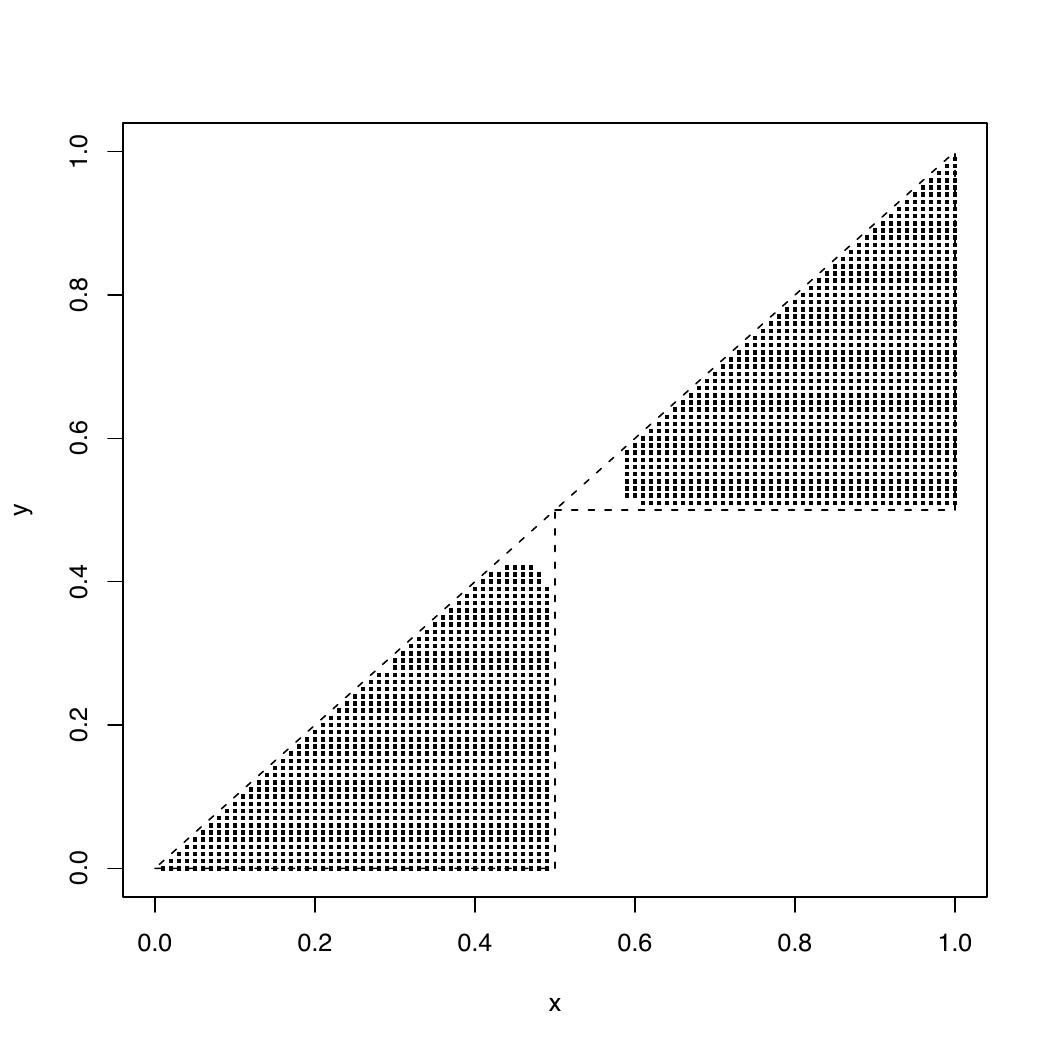}&
\includegraphics[width=0.23\textwidth ]{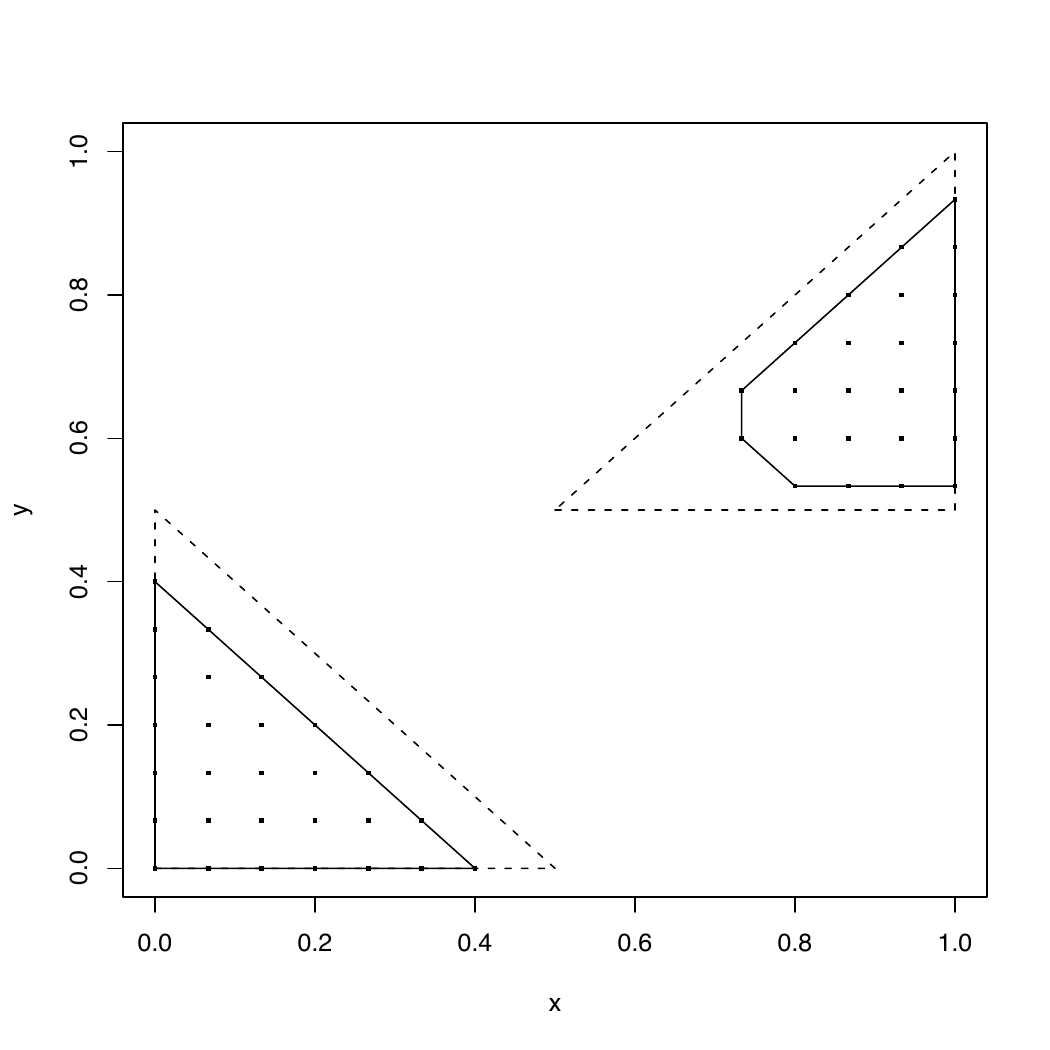}
\includegraphics[width=0.23\textwidth ]{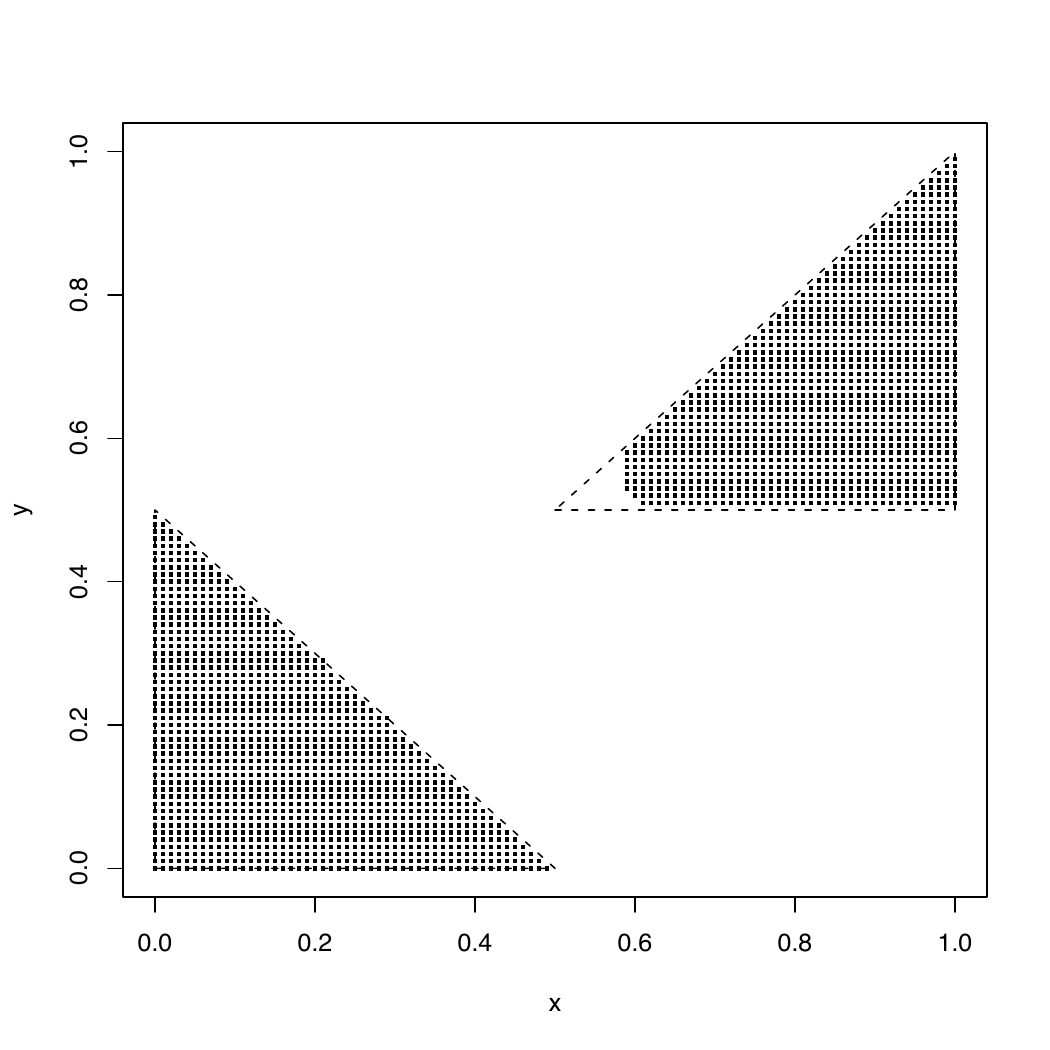}\\
(a)  & (b)
\end{tabular}
\caption{
(a) Approximations of $G$ by $G_{ a}$ with $a=(\alpha, \alpha^{\prime})$ and $\alpha =\alpha^{\prime} =15$ (left) and $\alpha =\alpha^{\prime} =100$ (right);
(b) Approximations of $G^{\prime}$ by $G^{\prime}_{a}$  with $\alpha =\alpha^{\prime} =15$ (left) and $\alpha =\alpha^{\prime} =100$ (right)}
 \label{Figure 2}
\end{center}
\end{figure}

\subsection{Approximation using the level sets}

To demonstrate the performance of the Eq. \eqref{Eq:2}, let us denote the level sets of $f$ by $G_f (\lambda)=\{{\bf x} \in [0, 1]^2: f ({\bf x}) \geq \lambda \}$, where $\lambda > 0$.  Assume
now $f=I_G$, and for each image $G$, let us denote by $G_{a}:=G_{f_{a}} (1/2)$ the  corresponding moment-recovered image with $f_{a}$ defined  in \eqref{Eq:2}.

\vspace{0.3cm}

\noindent
\textsc{\bf Example 2}.   Let $G_{\bf a}$ and  $G^{\prime}_{a}$ be the approximations of two regions $G=G_1 \cup G_2$ and $G^{\prime}=G^{\prime}_1 \cup G_2$, respectively,
representing the union of two right triangles located in different
orientations within $[0, 1]^2$.  See the plots in  Figure 2 where the boundaries of $G$ and $G^{\prime}$ are denoted by dashed lines. Namely, consider
\begin{align*}
&G_1=\{(x, y):  0 \leq x \leq 0.5,\, 0 \leq y\leq x  \} \,, G_2=\{(x, y):  0.5 \leq x \leq 1,\, 0.5 \leq y\leq x  \}  \; {\rm and} \\\notag
&G^{\prime}_1=\{(x, y) :  0 \leq x \leq 0.5,\, 0 \leq y\leq  0.5 - x  \}\, .
\end{align*}
At first, we took $\alpha =\alpha^{\prime} =15$ in  $G_{\bf a}$ and  in $G^{\prime}_{a}$ (the solid lines), and then, we took $\alpha =\alpha^{\prime}
=100$, see Figure 2(a)-2(b). Comparing the plots in Figure 2(a) and 2(b), we conclude  that the approximation around the point $(0.5, 0.5)$ is not as accurate as in  other locations.

The plots in Figure \ref{Figure 2} justify that the approximants  almost coincide with the target objects  as $\alpha =\alpha^{\prime} =100$. Let us  denote by $A  \Delta B$ the symmetric
difference between sets $A$ and $B$. We calculated the values of pseudo-metrics $\varepsilon_{ a}= \lambda (G_{a} \Delta G)$ induced by the Lebesgue measure
$\lambda$ on $[0, 1]^2$  for  $\alpha=\alpha^{\prime} \in \{15, 20, 25,  32, 50, 100\}$; see Table 2, where
$\varepsilon_{a}^{\prime}=\lambda(G_{a}^{\prime} \, \Delta G^{\prime})$.

\begin{table}
\centering
\caption{The  $\varepsilon_{\bf a}$- and $\varepsilon^{\prime}_{\bf a}$-errors}
{\begin{tabular}{|c|c|c|c|c|c|c|}
 \hline
  $\alpha=\alpha^{\prime}$ & 15 & 20 & 25  & 32 & 50 & 100 \\
 \hline
$\varepsilon_{\bf a}$ & 0.1233 & 0.1000 &  0.0724 & 0.06836 & 0.0452  & 0.02425  \\
  \hline
$\varepsilon_{\bf a}^{\prime}$ & 0.1012 & 0.07375 & 0.0644 & 0.04932 & 0.0326  & 0.01555   \\
  \hline
\end{tabular}}
\label{table2}
\end{table}
Also, one can compare the approximants in Figure 2(a) and 2(b) with the
approximations of right triangles derived in He  \cite{he2001}.  See also
Gustafsson {\it et al}  \cite{gustafsson2000}, where the approximations of
so-called quadrature domains (with a real analytic smooth boundaries)
via a formal exponential transform of the moment sequence are
derived. The approximations based on the exponential transform
mostly work for the cases when the dimension $d \leq 2$. Our
construction can be applied for recovering a support boundary of any high dimensional $M$-determinate function $f$ (not necessary to be $I_G$) by means of the  level sets $G_{f_a}$. 

\subsection{Recovery of $f=I_G$}
Now, let us apply Eq. \eqref{Eq:2} for approximation of a discontinuous function. Assume that $f$ is piece-wise constant with jumps on the boundary of a connected regular set $G$ with smooth boundary. For simplicity, assume that $f$ is the indicator function of $G$, i.e., let $f=I_G$.  To be able to derive the rate of approximation in this case, let us recall the beta kernel density functions $\beta_{\alpha} (\cdot, x)$ introduced in Remark 1, and rewrite them as follows:
$$
\beta_\alpha (t, x) = \frac{\Big (t^a (1-t)^b\big )^\alpha}{B(a\alpha +1, b \alpha + 1)},  \;\;\; a=\frac{[\alpha x]}{\alpha},\;\; b=1- \frac{[\alpha x]}{\alpha}, \;\;\; 0 \leq t, x \leq 1.
$$
Here $ B (\cdot, \cdot) $ is the beta function, $ a + b = 1 $, and $ \beta_\alpha (\cdot, x) $ approaches  the Dirac delta function $ \delta (\cdot -x) $ on $ [0,1] $ as $ \alpha \rightarrow \infty $.  We want to estimate the rate of convergence of this sequence to the delta function. For simplicity assume that  the shape parameters  of $\beta_{\alpha} (\cdot, x)$ density function  are equal to $\alpha x +1$ and $\alpha - \alpha x +1$, respectively. In this case, in  the previous expression for $\beta_\alpha (t, x)$ we will have $a=x$ and $b=1-x$. Denote  $D_{>}=\{u: | u-x | > d \}\cap [0, 1]$.

\vspace{0.05in}

\noindent
{\bf Proposition 2}. For each $x\in (0, 1)$, $ d < {\rm dist} (x, \partial[0,1])/2 $, and  sufficiently large $\alpha \in \mathbb {N}_+$, there exist $ B < 1 $ such that 
\begin{align}\label{Eq:3}
\int_{D_{>}}  \beta_{\alpha} (t, x) d t <  B^{\alpha} \sqrt{\frac{\alpha+2}{2 \pi x (1-x)  }} .
\end{align}

\medskip

\noindent
{\bf Proof.} Using Stirling's formula,  the beta function behaves asymptotically as
$$ 
B(u,v) \approx \sqrt{2 \pi} \frac{u^{u - \frac{1}{2}} v^{v- \frac{1}{2}}}{(u+v)^{u+v - \frac{1}{2}}} ,
$$
for large values of $ u $ and $ v $, and  
\begin{align}\label{Eq:4}
\beta_\alpha (t, x) \approx \sqrt{\frac{\alpha+2}{2 \pi x (1-x) }} \ g(t, x)^\alpha ,
\end{align}
where 
$$
g (t, x) = \Big(\frac{t}{a}\Big)^a \   \Big(\frac{1-t}{b} \Big)^b, \ a=x, \;\; b=1-x .
$$
(cf. with Bertin and Klutchnikoff \cite{bertin2011}). Taking logarithm of both sides of the above expression,
we have 
$$
\log \ g (t, x) = a \ \log(\frac{t}{a}) + b \ \log(\frac{1-t}{b}). 
$$
Since $ a+ b = 1 $ and $ \log t $ is a concave function,  
we will have
$$
\log g  (t, x) \leq \log (a \frac{t}{a} + b \frac{(1-t)}{b} )= \log 1 = 0.
$$
So, $ g (t, x) \leq 1 $ for all $ t $, or equivalently $$||g (\cdot, x)||_\infty \leq 1. $$

\medskip

Note that $ g (t, x) \geq 0 $ and is continuous.  If we use $ D_> $ to denote the exterior of a $ d$-neigborhood of $ x $ and $ D_< $ the interior of this neighborhood, i.e.
$$
\begin{array}{l}
D_> = \{ t: |t-x| > d \} \cap [0,1] \\
D_< = \{ x: |t-x| < d \} \cap [0,1],
\end{array}
$$
clearly
$$
\int_0^1 g (t, x) d t = \int_{D_<} g (t, x) d t + \int_{D_>} g (t, x) d t .
$$
By continuity and positivity of $ g (\cdot, x) $, both integrals on the right hand side are strictly greater than $ 0$. Also, since $ g(\cdot, x) $ has its peak at $ x $, for $ d $ ($ d$ depends on $ x $) chosen appropriately,
$$
\int_{D_>} g (t, x) d t < \frac{1}{2} \int_0^1 g (t, x) d t .
$$
The factor of $ \frac{1}{2} $ can be replaced by $ \delta < 1 $ by varying the choice of $ d $. By H\"{o}lder's inequality, since
$$
||g (\cdot, x)||_1 \leq ||g   (\cdot, x)||_\alpha \leq ||g  (\cdot, x) ||_\beta \leq ||g  (\cdot, x)||_\infty, \ \ \ 1 < \alpha \leq \beta < \infty, 
$$
we have 
$$
\Big (\int_{D_>} g^\alpha(t, x) dt \Big )^{\frac{1}{\alpha}}  <  \Big ( \int_0^1 g^\alpha (t, x) dt  \Big )^{\frac{1}{\alpha}} \leq ||g (\cdot, x)||_\infty \leq 1. 
$$
Since the left-hand side is strictly less than $ 1 $, there exist $ B < 1 $ so that $ (\int_{D_>} g^\alpha (t, x) d t)^{\frac{1}{\alpha}} < B  $, and
$$
\int_{D_>} g^\alpha (t, x) d t < B^\alpha .
$$
 Hence, for every $ x \in (0,1) $ and $d$ chosen appropriately (depending on $ x $), we obtain 
\begin{align}
\label{Eq:5}
\int_{D_>} g^\alpha (t, x) d t < B^\alpha. 
\end{align}
Note that if $ d = {\rm dist} (x, \partial[0,1])/2 $, $ B = 1- d $ satisfies \eqref{Eq:5}. 
Now combining \eqref{Eq:4} and  \eqref{Eq:5} yields  \eqref{Eq:3}.
\hfill{$\blacksquare${}}

\medskip

\noindent
{\bf Proposition 3}.  {\it Let $f=I_{G}$ for a closed convex subset $G \subseteq [0,1]^2$.  Then for each   ${\bf x} =(x, y)\in (0, 1)^2$, $ d \leq {\rm dist} ({\bf x}, \partial G)/2  $, and  sufficiently large values of $\alpha, \alpha^{\prime} \in \mathbb {N}_+$, there exist $ B < 1 $ such that}
\begin{align}
\label{Eq:6}
| f_{a} ({\bf x}) - f ({\bf x}) | \leq \frac{ B^{\alpha +\alpha^{\prime}}}{2 \pi} 
  \sqrt{\frac{(\alpha+2)\,(\alpha^{\prime}+2)}{x\,(1-x) y\,(1-y)}}
  \end{align}
{\it and $f_{a} \longrightarrow_{L_1}  f $  as  $a\to\infty$, with}
\begin{align}
d_1 ( f_{a}, f)&\leq \frac{4}{\pi \; 2^{\alpha + \alpha^{\prime}}} 
 \frac{1}{\sqrt{(\alpha +1/2) (\alpha'+ 1/2)}}
  \end{align}
  
  \noindent
{\bf Proof.} Since  $f=I_G$, after changing the order of summation and integration in  \eqref{Eq:2}, we obtain, for each ${\bf x}=(x, y) \in (0, 1)^2$ 
\begin{align}\label{Eq:7}
&f_{a} ({\bf x}) - f ({\bf x})  =  \int_{[0,1]^2} \,   \beta_\alpha (t, x) \, \beta_{\alpha^\prime} (s,  y) [I_G (t, s) - I_G (x, y)] \, d t \, d s\\\notag
&=  \int_{[0,1]^2} K_a ({\bf x}, {\bf t})  [I_G ({\bf t}) - I_G ({\bf x})] d {\bf t} = I_{G^c} ({\bf x})   \int_{G} K_a ({\bf x}, {\bf t}) d {\bf t}  - I_{G} ({\bf x})   \int_{G^c} K_a ({\bf x}, {\bf t}) d {\bf t},
\end{align}
since the difference
\begin{equation}\label{Eq:diff}
I_G ({\bf t}) - I_G ({\bf x})= 
\begin{cases}
I_{G} ({\bf t}), \quad \;\;\;\;   {\bf x} \in G^c= [0,1]^2\setminus G,\\
 - I_{G^c} ({\bf t}), \quad {\bf x} \in G.\\
\end{cases}
\end{equation}
Here $K_a ({\bf x}, {\bf t}):= \beta_\alpha (t, x) \, \beta_{\alpha^\prime} (s,  y)$. Let us estimate the first term in the last equation of \eqref{Eq:7}. For each ${\bf x}\in G^c$, there exists a small disc with radius $d$ such that $G\subset D_{>}=\{{\bf u}:  || {\bf u} - {\bf x} || > d\}$. Hence, for each ${\bf x} \in G^c$, we can write
\begin{equation}\label{Eq:9}
\int_{G} K_a ({\bf x}, {\bf t}) d {\bf t} \leq \int_{D_>} K_a ({\bf x}, {\bf t}) d {\bf t} \leq  \int_{D^x_{(>, \frac{d}{\sqrt 2})}}  \beta_\alpha (t, x) dt \int_{D^y_{(>, \frac{d}{\sqrt 2})}}  \beta_\alpha (s, y) ds,
\end{equation}
where $D^x_{(>, \frac{d}{\sqrt 2})}$ is the exterior domain in x-direction of radius $\frac{d}{\sqrt 2}$, and similarly $D^y_{(>, \frac{d}{\sqrt 2})}$ is the exterior domain in y-direction of radius $\frac{d}{\sqrt 2}$. Application of Proposition 2 twice in the right-hand side of  \eqref{Eq:9}  yields \eqref{Eq:6}. 

\medskip

To determine the $L_1 $ rate of convergence $ f_a $ to $f $, recall that $ B $ in \eqref{Eq:6} depends on the distance of $ {\bf x} $ to the boundary of $ G $. Considering the projection of $ G \subset [0,1]^2 $ onto each of the coordinate axes, if  $ x \in [0, \frac{1}{2}] $, we may choose $ d = x $ and if $ x \in [\frac{1}{2}, 1] $, then $ d = {1-x} $. Hence by symmetry, the exterior integral of the beta function in \eqref{Eq:6} will satisfy 

\begin{align}\label{Eq:12}
\int_{\frac{1}{2}}^{1}  \frac{(1-x)^{\alpha- \frac{1}{2}}}{\sqrt{x}} dx &\leq  \sqrt{2} \int_{\frac{1}{2}}^{1} (1-x)^{\alpha - \frac{1}{2}} dx \\
& \leq \frac{1}{2^{\alpha}\, (\alpha + \frac{1}{2})},
\end{align}
if $ \alpha > 1 $. Integrating \eqref{Eq:6} on $[0,1]^2$, since $ G $ is convex, we may decompose $ d({\bf x}, \partial G) $ into its components. Using the symmetry argument leading to \eqref{Eq:12}
\begin{align}
\int_{[0,1]^2} |f_a ({\bf x}) - f({\bf x}) |d{\bf x} & \leq \frac{2}{\pi} \sqrt{(\alpha+2)(\alpha^{\prime} +2)} \int_0^1 \int_0^1 \frac{d(x, \partial[0,1])^\alpha d(y, \partial [0,1])^{\alpha^{\prime}} }{\sqrt{x(1-x)y(1-y)}} dx dy \\
&= \frac{8}{\pi} \sqrt{(\alpha+2)(\alpha'+2)} \int_{\frac{1}{2}}^{1} \int_{\frac{1}{2}}^{1}  \frac{(1-x)^{\alpha - \frac{1}{2}}}{\sqrt{x}}  \frac{(1-y)^{\alpha^{\prime} - \frac{1}{2}}}{\sqrt{y}} dx dy \\
& \leq  \frac{4}{\pi}  \frac{(1/2)^{\alpha + \alpha^{\prime}}}{\sqrt{(\alpha + \frac{1}{2})(\alpha' + \frac{1}{2} )}} \ \ {\rm as} \ a= (\alpha, \alpha') \rightarrow \infty .
\end{align}
\hfill{$\blacksquare${}}

\begin{figure}
\begin{center}
\begin{tabular}{ccc}
 \includegraphics[width=0.35\textwidth ]{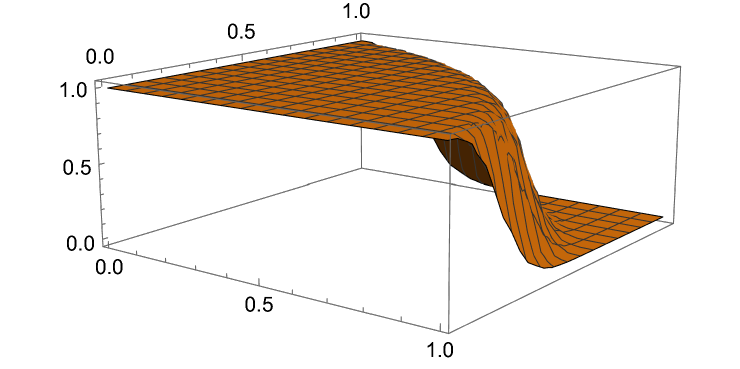} &
 \includegraphics[width=0.4\textwidth ]{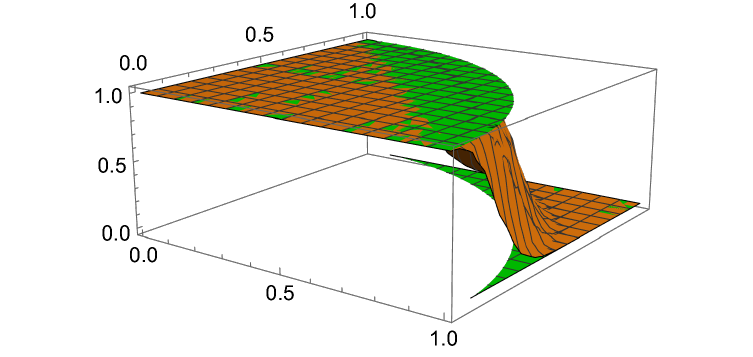}&\\
(a)  & (b)
\end{tabular}
\caption{\small (a) Approximation $f_{ a}$ of $f=I_{G_1}$, where $G_1$ is the quarter of a unit disc with the center at  origin; (b) $f=I_{G_1}$  (green)
and its approximant $f_{a}$, with $\alpha=\alpha^\prime =75$}
\label{Figure 3}
\end{center}
\end{figure}


\noindent
{\bf Corollary 3}.  {\it Let $f=c_1 I_{G_1} + c_2 I_{G_2}$ be a linear combination of two indicator functions of disjoint sets $G_1$ and $G_2$.  Then $f_{a}  \longrightarrow_{L_1} f$   as $a\to\infty$.}

\vspace{0.05in}

\noindent
{\bf Example 3}. Consider the indicator function $f=I_{G_1}$, where $G_1$ represents  a quarter of a unit disc $x^2 +y^2\leq 1$. Application of (2) with geometric moments
$$
\mu_{G_1} (k, j) = \frac{1}{2\, (j+1)}\, B\Big( \frac{k+1}{2}, \frac{j+3}{2}\Big), \;\;\;  (k, j) \in \mathbb{N} _{a}, 
$$
where $B(u, v)$ denotes the Beta function with the shape parameters $u$ and $v$,
provides the approximation $f_{a}$ of $f$ as $a\to\infty$. See the plots  in Figure 3. In a similar  way the approximation  of $f=I_{G_2}$ is constructed,  when $G_2$ is a half of a  disc  with a center at $(0.5, 0)$ and the  radius equal to $0.5$. See the plot of  approximation  of $f$ displayed  in Figure 4 (a). Here we set  $\alpha=\alpha^\prime =500$.

\begin{figure}[!ht]
\begin{center}
\begin{tabular}{ccc}
 \includegraphics[width=0.35\textwidth ]{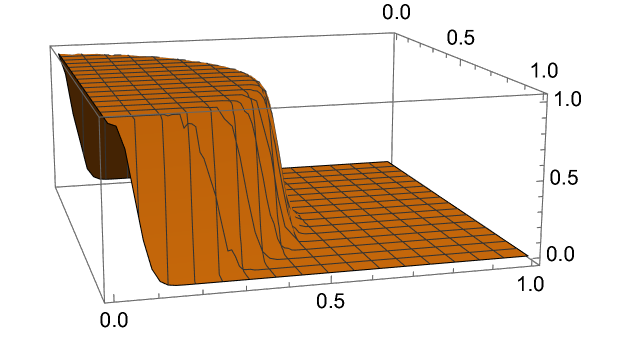} &
 \includegraphics[width=0.35\textwidth ]{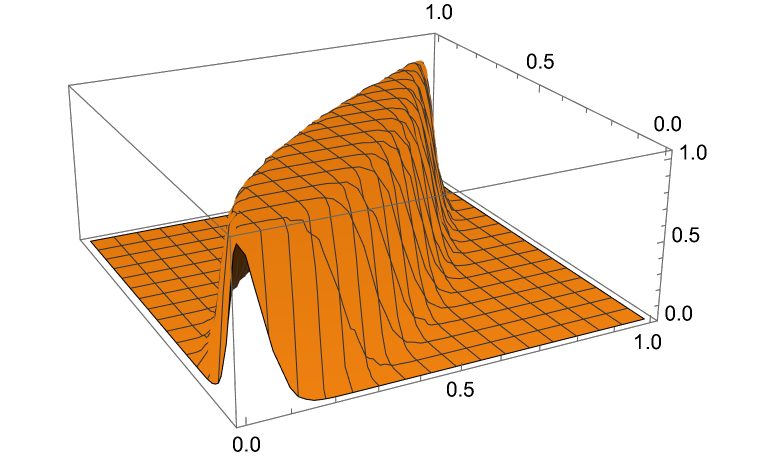}&\\
(a)  & (b)
\end{tabular}
\caption{\small (a) Approximation $f_{ a}$ of $f=I_{G_2}$, where $G_2$ is the half of a  disc with the center at  (0.5, 0) and   radius 0.5 and $\alpha=\alpha^\prime =500$; (b) approximation of $f=I_{G_3}$ where  $G_3=\{(x, y): 0\leq x \leq 1, x^2 \leq y \leq \sqrt x\}$ with $\alpha=\alpha^\prime =200$}
\label{Figure 4}
\end{center}
\end{figure}

\noindent
{\bf Remark 2.} To improve the performance of our approximation, let us consider  the following modification  of $f_a$ (cf. with Corollary 1 from Mnatsakanov and Garai \cite{mnatsakanov2017}, where similar modification has been used in univariate case). Namely, let us denote by $\widetilde a$ the  two-dimensional vector of parameters $\widetilde a=(\alpha/c, \alpha^{\prime}/c)$, and introduce the following  linear combination of $f_{\widetilde a}$ and $f_{a}$:
\begin{align}
{\widetilde f}_{a}:= \frac{1}{1-c} f_{\widetilde a} -\frac{c}{1-c} f_{a},
\end{align}
with some constant  $0 < c < 1$. In Table 3, we can see the improved records of corresponding    $L_1$-rates for ${\widetilde f}_{a}$, when $c=0.5$,  while using the same order of geometric moments of unknown set $G_1$ if compared to  $f_{a}$. Here $G_1$ represents  a quarter of a unit disc $x^2 +y^2\leq 1$ with a center at origin (introduced in Example 3).

\begin{table}
\caption{\small   Example 3. Approximation errors   $d_1 (f_{a}, f)$  and $d_1 (\widetilde {f}_{a}, f)$ with $c=0.5$.}
{\small
\begin{center}
\begin{tabular}{|l|rrrrrrr|}
\hline
\;\;\;\;\;\;\; $\alpha=\alpha^{\prime}$ & 20 & 40 & 60 & 80 & 100 & 120 & 140  \\
\hline
$d_1 (f_{a}, f)\times 10^3$& 51.92& 25.52& 16.91& 12.64& 10.09& 8.396& 7.19\\
\hline
 \hline
\;\;\;\;\;\;\; $\alpha=\alpha^{\prime}$ & 10 & 20 & 30 & 40 & 50 & 60 & 70  \\
\hline
$d_1 (\widetilde {f}_{a}, f)\times 10^3$& 2.50& 0.63& 0.28& 0.16& 0.10& 0.069& 0.051\\
\hline
\end{tabular}
\end{center}
}
\end{table}

\vspace{0.05in}

\noindent
\text {\bf  Example 4}. Consider the indicator function $f=I_{G_3}$  of the set $ G_3 $  bounded by curves $y=x^2$ and $y=\sqrt x$, i.e.,  $G_3=\{(x, y): 0\leq x \leq 1, x^2 \leq y \leq \sqrt x\}$.
The plot of approximation of $f=I_{G_3}$  is presented in Figure 4 (b), when $\alpha=\alpha^\prime =200$.  In this case the sequence of  geometric moments of $G_3$ is specified as follows:
\begin{equation*}
m_{G_3} (k, j)= \frac{3}{(2k + j +3)\,(k + 2j +3)}, \;\;\;\; (k, j)\in \mathbb{N}_{a}.
\end{equation*}



\section{Recovering the discrete distributions with finite support}

Assume that $f$ is a discrete probability density function with support  $supp\{f\}=\{x_k, \ k=0, \dots, N \}$. Assume  $x_k\in [0, 1]$. Our goal is to recover 
$f(x_k)=P\{x=x_k\}:=p_k $ or corresponding cumulative distribution function $F (x)=\sum_{k=0}^N p_k I_{[0,  x)} (x_k)$ given the sequence of its moments
\begin{equation}
\label{Eq:mom}
m (j, dF) = \sum_{k=0}^N \,  x_k^j  p_k =\int_0^1 t^j d F (t),\;\;\; j=0, 1, \dots, \alpha, 
\end{equation}
known up to integer order $\alpha\in \mathbb{N}_+$.  Consider the following approximation of $F$:
\begin{equation}
F_{\alpha} (x) =(\mathcal{K}^{-1}_{\alpha} m (\cdot, dF) ) (x), \;\;\; x \in [0, 1].
\end{equation}
Here, for each $0 \leqslant  x \leqslant 1$ and $\alpha \in \mathbb{N}_+$
\begin{equation}
\label{Eq:cdf}
\left(\mathcal{K}^{-1}_{\alpha} m (\cdot, dF)\right) (x) = \sum^{[\alpha x]}_{k=0} \sum^{\alpha }_{j=k} \,
\binom{\alpha }{j} \binom{j}{k}\,  (-1)^{j-k} \, m (j, dF).
\end{equation}

\noindent
{\bf Proposition 4}. {\it For each $x \in [0, 1]$,  we have $F_{\alpha} (x) \to F(x)$ as $\alpha \to\infty$}.

\vspace{0.04in}

\noindent
{\bf Proof.} Indeed, after substituting  \eqref{Eq:mom} into \eqref{Eq:cdf},  we have
\begin{multline}
\label{Eq:conv}
F_{\alpha} (x) =(\mathcal{K}^{-1}_{\alpha}\, m (\cdot, d F)) (x) = \int_0^1 \, \sum^{[\alpha x]}_{k=0}  
\sum^{\alpha }_{j=k} \, \binom{\alpha }{j} \binom{j}{k}\,     \,
t^k \, \big(-t \big)^{j - k}
\, d F(t)\\
= \int_0^1 \sum^{[\alpha {x}]}_{k=0}\binom{\alpha }{k} \, t^k 
\sum^{\alpha -k}_{m=0} \,  \binom{\alpha  - k}{m} \big(-t\big)^m
\, d F(t)
= \int_0^1\sum^{[\alpha {x}]}_{k=0}\binom{\alpha }{k} \,t^k   (1- t)^{\alpha-k}  d F(t) \\
=\int_0^1 B_{\alpha } \big (t, x \big)\, d F(t) \shoveleft{\to \int_0^1 1_{[0,x)} (t)\, d F(t) = \sum_{i=0}^N 1_{[0,x)} (x_i)\,  p_i= F (x), \quad \text{as $\alpha \to \infty$}.}
\end{multline}
The convergence in \eqref{Eq:conv}  is justified by using simple algebra and the properties of  the binomial probabilities 
\begin{equation}\label{Eq:Bin}
B_{\alpha} (u, v) = \sum^{[\alpha v]}_{k=0} \binom{\alpha}{k}   u^k \, \big (1 - u\big)^{\alpha  - k}\to
\begin{cases}
1, \quad u < v\\
0, \quad u > v
\end{cases},\; \alpha \to \infty.
\end{equation}
\hfill{$\blacksquare${}}

To recover the value of $f (x_j)$ one can take $f_\alpha (x_j)= F_{\alpha} (x_{j+1})- F_{\alpha} (x_{j})$.

\noindent
 Similar statement is true for multidimensional case. Indeed, assume that $ supp \{f\}$ represents a finite (may be countable) subset of $[0, T]^2$ for some $0 < T < \infty$. Without loss of generality, assume  $T=1$, and denote $supp \{f\}:=\{(x_i, y_j), 0 \leq i, j \leq N\}$, and define the following approximation of 
$$
F (x, y)=\sum_{i=0}^N \sum_{j=0}^N P(X=x_i,Y= y_j) 1_{[0,x)} (x_i) 1_{[0,y)} (y_j)=\sum_{i=0}^N \sum_{j=0}^N p_{i,j} 1_{[0,x)} (x_i) 1_{[0,y)} (y_j)
$$
constructed as follows
\begin{equation}\label{Eq:2.0}
F_{a} ({\bf x}) :=\left(\mathcal{K}^{-1}_{a}  \nu \right)({\bf x}) = \sum^{[\alpha {x}]}_{k=0} \sum^{[\alpha^\prime {y}]}_{l=0} \sum^{\alpha }_{j=k}  \sum^{\alpha^\prime }_{n=l} 
\binom{\alpha }{j} \binom{j}{k} \binom{\alpha^\prime}{n} \binom{n}{l}(-1)^{j+n- k-l} m  (j, n, d F),
\end{equation}
for each ${\bf x}= (x, y)\in [0,1]^2$. Here, by $\nu$ we denote the sequence of so-called product moments, $\nu= \{m ({j, k}, d F), (j, k)\in \mathbb {N}_a\}$, where 
\begin{equation}\label{Eq:2.1}
m ({j, k}, d F) = \int_0^1\int_0^1 t^j\, s^k d  F(t, s), \; \;  j, k \in \mathbb{N}= \{0, 1, \dots \}, \;\; m  (0,0)=1.
\end{equation}

\noindent
{\bf Proposition 5}.
{\it For each ${\bf x} \in [0, 1]^2$,  we have $F_{\alpha} ({\bf x}) \to F({\bf x})$ as $a \to\infty$}.

\vspace{0.04in}

 \noindent
{\bf Proof.} Indeed, consider the sequence of  operators $\mathcal{K}^{-1}_{a}$: 
\begin{equation}\label{Eq:2.2}
\left(\mathcal{K}^{-1}_{a}\,  \nu \right) (x, y) = \sum^{[\alpha {x}]}_{k=0} \sum^{[\alpha^\prime {y}]}_{l=0} \sum^{\alpha }_{j=k}  \sum^{\alpha^\prime }_{n=l} \,
\binom{\alpha }{j} \binom{j}{k}\, \binom{\alpha^\prime}{n} \binom{n}{l}\, (-1)^{j+n- k-l}  \,m  (j, n, d F), 
\end{equation}
where $0 \leq  x, y \leq 1\, , \, \,  a=(\alpha , \alpha^\prime)$ with $\alpha,  \alpha^\prime  \in \mathbb{N}_+$. The sequence of operators $\mathcal{K}^{-1}_{a}\,  \nu$ is determined by  the values of ordinary moments \eqref{Eq:2.1}.  
 In the sequel, we write  $a \to \infty$  to  mean  that $\alpha \to\infty$ and $ \alpha^\prime \to \infty$, and $\mathbb{N}_a:=\mathbb{N}_{\alpha}\times \mathbb{N}_{\alpha^{\prime}}$ with $\mathbb{N}_{\alpha}= \{0, 1, \dots, \alpha\}$. 
From Mnatsakanov (\cite{mnatsakanov2011}, Thm.1) we know that 
$\mathcal{K}^{-1}_{a} \, \mathcal{K} F  \longrightarrow_w F, \, \, {\rm as} \, \, a \to\infty$, where ``$\longrightarrow_w$" denotes convergence
of cdfs in a weak sense, i.e., the convergence at each continuity point of the limiting cdf $F$.
In a similar way, as on the second line of  \eqref{Eq:conv}, we can rewrite $F_{a} (x, y)$ and  prove the statement of Proposition 5 by  using the Lebesgue dominated convergence theorem:
\begin{multline*}
F_{a} (x, y) = \int_0^1\int_0^1  \sum^{[\alpha {x}]}_{k=0}  
\sum^{\alpha }_{j=k} \, \binom{\alpha }{j} \binom{j}{k}\,     \,
t^k \, (-{t})^{j - k}\sum^{[\alpha^\prime {y}]}_{l=0}\sum_{n=l}^{\alpha^\prime}\binom{\alpha^\prime}{n} \binom{n}{l}\, 
\,s^l \, (-{s})^{n- l} d F(t, s)\\
= \int_0^1\int_0^1 \sum^{[\alpha {x}]}_{k=0}\binom{\alpha }{k} \,t^k 
\sum^{\alpha -k}_{m=0} \,  \binom{\alpha  - k}{m}  (- t)^m \, \sum^{[\alpha^\prime {y}]}_{l=0}\binom{\alpha^\prime }{l} \,s^l 
\sum^{\alpha^\prime -l}_{m^\prime=0} \,  \binom{\alpha^\prime  - l}{m^\prime}  (- s)^{m^\prime}
\, d F(t, s)\\
= \int_0^1\int_0^1 B_{\alpha }  ({t}, {x})\, B_{\alpha }  ({s}, {y})d F(t, s) \\
\shoveleft \to   \int_0^1\int_0^1 1_{[0,x)} (t) 1_{[0,y)} (s)  d F(t, s)
=\sum_{k=0}^N \sum_{j=0}^{N}1_{[0,x)} (x_k) 1_{[0,y)} (y_j)\,  p_{k, j}= F (x, y),
\end{multline*}
as $a \to \infty$.\hfill{$\blacksquare${}}

\vspace{0.25in}

\centerline{ \bf Compliance with Ethical Standards}


\noindent
{\bf Funding.} The authors did not receive support from any organization for the submitted work.


\noindent
{\bf Disclosure of potential conflicts of interest.} The authors declares no competing interests.


\noindent
{\bf Data Availability Statement.} No data has been used for this study.


\begin{thebibliography}{10}

\bibitem{Akheizer2020}
N.~I. Akheizer, \emph{The classical moment problem}, Wiley, New York, 1970.

\bibitem{bertin2011}
K.~Bertin and N.~Klutchnikoff, \emph{Minimax properties of beta kernel estimators}, Journal of Statistical Planning and Inference \textbf{141} (2011), 2287--2297.

\bibitem{bouezmarni2003}
T.~Bouezmarni and J.-M. Rolin, \emph{Consistency of the beta kernel density function estimator}, Canadian Journal of Statistics \textbf{31} (2003), 89--98.

\bibitem{chen1999}
S.~X. Chen, \emph{Beta kernel estimators for density functions}, Computational Statistics \& Data Analysis \textbf{31} (1999), 131--155.

\bibitem{Choi2020}
Hayoung Choi, Victor Ginting, Farhad Jafari, and Robert Mnatsakanov, \emph{Modified radon transform inversion using moments}, Journal of Inverse and Ill-posed Problems \textbf{28} (2020), no.~1, 1--15.

\bibitem{Cuyt2005}
Annie Cuyt, Gene Golub, Peyman Milanfar, and Brigitte Verdonk, \emph{Multidimensional integral inversion, with applications in shape reconstruction}, SIAM Journal on Scientific Computing \textbf{27} (2005), no.~3, 1058--1070.

\bibitem{goldenshluger2004}
A.~Goldenshluger and V.~Spokoiny, \emph{On the shape-from-moments problem and recovering edges from noisy radon data}, Probability Theory and Related Fields \textbf{128} (2004), 123--140.

\bibitem{gustafsson2000}
B.~Gustafsson, C.~He, P.~Milanfar, and M.~Putinar, \emph{Reconstructing planar domains from their moments}, Inverse Problems \textbf{16} (2000), 1053--1070.

\bibitem{he2001}
C.~He, \emph{Moment problems and operator theory}, Ph.D. thesis, University of California, 2001.

\bibitem{henrion2023}
D.~Henrion, M.~Korda, and J.-B. Lasserre, \emph{Polynomial argmin for recovery and approximation of multivariate discontinuous functions},  (2023), arXiv:2302.06945.

\bibitem{mnatsakanov2011}
R.M. Mnatsakanov, \emph{Moment-recovered approximations of multivariate distributions: The laplace transform inversion}, Statistics and Probability Letters \textbf{81} (2011), 1--7.

\bibitem{mnatsakanov2022}
R.M. Mnatsakanov, H.~Albrecher, and S.~Loisel, \emph{Approximations of copulas via transformed moments}, Methodology and Computing in Applied Probability \textbf{24} (2022), 113557.

\bibitem{mnatsakanov2017}
R.M. Mnatsakanov and B.~Garai, \emph{On the moment-recovered approximations of regression and derivative functions with applications}, Journal of Computational and Applied Mathematics \textbf{315} (2017), 17--31.

\bibitem{Shohat}
J.A. Shohat and J.D. Tamarkin, \emph{The problem of moments}, American Mathematical Society, New York, 1943.

\end{thebibliography}



\providecommand{\bysame}{\leavevmode\hbox to3em{\hrulefill}\thinspace}
\providecommand{\MR}{\relax\ifhmode\unskip\space\fi MR }
\providecommand{\MRhref}[2]{%
  \href{http://www.ams.org/mathscinet-getitem?mr=#1}{#2}
}
\providecommand{\href}[2]{#2}

\end{document}